\documentclass[a4paper]{report}
\usepackage[utf8]{inputenc}
\usepackage[T1]{fontenc}
\usepackage{RJournal}
\usepackage{amsmath,amssymb,array}
\usepackage{booktabs}
\usepackage{lscape}
\usepackage{multirow}
\usepackage{caption}
\usepackage{subcaption}

\begin{document}

\sectionhead{Contributed research article}
\volume{XX}
\volnumber{YY}
\year{20ZZ}
\month{AAAA}

\begin{article}

\title{\pkg{nlstac}: Non-Gradient Separable Nonlinear Least Squares Fitting}
\author{by J. A. F. Torvisco, R. Benítez, M. R. Arias and J. Cabello~Sánchez}

\maketitle

\abstract{
A new package for nonlinear least squares fitting is introduced in this paper. This package implements a recently developed algorithm that, for certain types of nonlinear curve fitting, reduces the number of nonlinear parameters to be fitted. One notable feature of this method is the absence of initialization which is typically necessary for nonlinear fitting gradient-based algorithms. Instead, just some bounds for the nonlinear parameters are required. Even though convergence for this method is guaranteed for exponential decay using the max-norm, the algorithm exhibits remarkable robustness, and its use has been extended to a wide range of functions using the Euclidean norm. Furthermore, this data-fitting package can also serve as a valuable resource for providing accurate initial parameters to other algorithms that rely on them.}

\section{Introduction}
Experimental data often exhibits non-linear patterns. As such, researchers in applied science often have to try to fit these data with non-linear models which can be challenging to fit. In this paper, we introduce the \CRANpkg{nlstac} package \citep{package_nlstac}, which implements the TAC algorithm for solving separable nonlinear regression problems, among others. Unlike other solvers, it does not require initialization values. Throughout the paper, we emphasize the potential synergistic usage of \pkg{nlstac} alongside other commonly used solvers, as it can provide reliable initialization values for them. The syntax of \pkg{nlstac} follows a similar structure to the solvers in the \CRANpkg{minpack.lm} package \citep{package_minpacklm}, making it familiar to researchers experienced with those solvers.

The motivation behind developing the \pkg{nlstac} package stems from an approximation problem involving time series data. Specifically, we were working with data corresponding to measurements obtained from a thermometer reaching thermal equilibrium with the surrounding medium, particularly oceanic water. In accord with Newton's law of cooling, the temporal evolution of these data exhibits an exponential pattern described by the expression:
\begin{equation}\label{eq:pp}
a_1 e^{-k_1t} + a_2, \; \text{where } a_1,a_2,k_1 \in \mathbb{R}, \; k_1>0,
\end{equation}
where $t$ represents the time variable.

Fitting data with the exponential pattern described by equation (\ref{eq:pp}) is a nonlinear optimization problem. One challenge we encountered with widely used algorithms for such problems was the need to initialize the parameters in (\ref{eq:pp}), as the solutions often strongly depended on the chosen initial values ---a bad choice of initial values could lead to a sub-optimal local minimum or even make the algorithm not to converge at all. The TAC algorithm, around which the present package is built, is presented in~\citet{tac} and it overcomes this issue by eliminating the requirement for parameter initialization. It only needs to specify a broad interval in which to search for the nonlinear parameters. As we worked with the TAC algorithm, its robustness became increasingly evident---robustness in the sense of stability of the algorithm in relation with noisy data and the convergence for a great variety of problems. In our opinion, this advantage, along with the lack of initialization, outweighs the need to specify the exact pattern to be used.

While the convergence of TAC is proven using the max norm, as shown in~\citet{tac}, we employ the Euclidean norm in the \pkg{nlstac} package and consider more general patterns beyond equation \eqref{eq:pp}. This extension of TAC beyond its proven convergence conditions is supported by its reliable performance, as mentioned earlier.

We acknowledge the widespread use of other algorithms for nonlinear fitting, such as Gauss-Newton or Levenberg-Marquardt, with the former being the default choice for the \code{nls} function in the \pkg{stats} package \citep{R}. Hence, in the present paper, we aimed to showcase the similarities and differences between \pkg{nlstac} and the \code{nls} fit, not as a competition, but as a demonstration of how well these two algorithms can work in synergy. Researchers sometimes encounter difficulties in finding suitable initialization values for \code{nls} to achieve convergence. In this regard, since \pkg{nlstac} does not require users to specify good starting values to converge, the resulting estimates can be used to provide good starting values to \code{nls} or any other initialization-dependent algorithm.

\subsection{Nonlinear regression, or more broadly, a problem of approximation}
Nonlinear regression or nonlinear fitting is a standard procedure commonly used in a wide range of scientific fields. Typically, we start with a dataset of $q$ observations of $n$ regressors (or predictor variables) and a response variable, namely $\{(\mathbf{x}_i, y_i),\, i = 1,\ldots,q\}$, and a mathematical expression relating the regressors  and the response variable. This mathematical expression may present a nonlinear dependency on several parameters. For instance, the aforementioned general expression is usually given by
\begin{equation}\label{eq:fullnonlinear}
    y = f(\mathbf{x};\boldsymbol\theta) + \varepsilon(\mathbf{x};\boldsymbol\theta),
\end{equation}
being $f:\mathbb{R}^n\times\mathbb{R}^p\rightarrow \mathbb{R}$ a function,  $\varepsilon(\mathbf{x};\cdot)$, independent and identically distributed random variables following a spherical normal distribution and  $\boldsymbol\theta = (\theta_1,\ldots,\theta_p)$ being the vector of parameters.
Thus, the problem reduces to finding an estimate of the parameter vector $\boldsymbol\theta^*$ such that some cost function is minimized. A usual choice for the cost function is the well-known least-squares cost function so that we are solving an optimization problem. Namely, finding $\boldsymbol\theta^* \in \mathbb{R}^p$ such that
\begin{equation}\label{eq:optimization}
    g(\boldsymbol\theta^*)=\min_{\boldsymbol{\theta}} g(\boldsymbol{\theta}), \;\;\; g(\boldsymbol{\theta})=\sum_{i = 1}^q\left(y_i-f(\mathbf{x}_i;\boldsymbol\theta)\right)^2.
\end{equation}

Please observe it is assumed that only $y_i$'s are observed with error, whereas $x_i$'s are measured exactly. The possibility of generalizing this to measurement-error models is not implemented in \pkg{nlstac}. More information about Least Squares Problems can be found in, for example,~\citet{bjorck_numerical_methods} or~\citet{nocedal_numerical_optization}.

The above problem can be described in Mathematical Analysis as an approximation problem. This kind of problem is determined by three elements: a set $A$, which is the object to be approximated; a family of functions $\mathcal{F}$, whose elements are known as approximants; and finally an approximation criterion---a procedure to measure how close to $A$ each element of $\mathcal{F}$ lies. In an approximation problem, a canonical question arises: does exist an element $f \in \mathcal{F}$ which is closest to $A$ ---attending to the approximation criterion--- than every other element in $\mathcal{F}$? When the answer to this question happens to be affirmative, a method to locate one of such elements is needed and this is what \pkg{nlstac} is designed for.

Let us identify those elements in the problem described above. The element $A$ to be approximated is the dataset of observations  $\{(\mathbf{x}_i, y_i),\, i = 1,\ldots,q\}$, which is a subset of $\mathbb{R}^n\times \mathbb{R}$. The family of approximants, $\mathcal{F}$, is given by the following $p$-parametric family
\begin{equation*}
\mathcal{F} = \left\{f_{\boldsymbol{\theta}}:\mathbb{R}^n\longmapsto \mathbb{R}\,|\,
f_{\boldsymbol{\theta}}(\mathbf{x}) = f(\mathbf{x};\boldsymbol{\theta}),\,\boldsymbol{\theta}\in\mathbb{R}^p\right\},
\end{equation*}
being $f$ the mathematical expression relating the regressors and the response variable given in \eqref{eq:fullnonlinear}. Finally, the approximation criterion corresponds to the function $g$ in \eqref{eq:optimization}. From now on, we will refer to one or the other definition depending on which one is more clear within the context.

Algorithms for finding the best set of parameters $\boldsymbol\theta^*$ of \eqref{eq:optimization} can be divided into local solvers and global solvers, depending on whether they are designed to find a local or a global minimum of the cost function, respectively. In general, local solvers are, under certain conditions, fast and accurate. They are usually based on some sort of gradient descent algorithm and are iterative in nature. That is, they start at a given initial guess for $\theta$ and, at each iteration (hopefully) they find a better approximation of the minimum. Under some assumptions (e.g. convexity), the local and global minima may coincide. However, in many cases, we will find that there are many different local minima and, consequently, the solution given by those algorithms may depend on the initial guess. Therefore, a bad initial guess could land us in a sub-optimal local minimum and there are even cases for which that initial conditions may cause the algorithm to not converge. Some local algorithms are the steepest descent method, incremental steepest descent method, Newton’s method, Quasi-Newton methods, Newton’s methods with Hessian modification, BFGS algorithm, Gauss-Newton method, or Levenberg-Marquardt method. More information about this and other methods can be found in, for example,~\citet{nocedal_numerical_optization},~\citet{arora_optimization_algorithms} or~\citet{rhinehart_nonlinear_regression_modeling}.

On the other hand, global solvers do not depend that heavily on an initial condition, but they require an interval or area in which to start looking for the minimum. Some global algorithms, like grid-search, are known to converge in any case, but they scale very poorly, such that the computational time grows exponentially with the number of parameters to be found. There are also heuristic solvers, which are not guaranteed to converge to global minimum, but they give a reasonable approximation in cases where other algorithms either take too long or do not converge at all. Some examples of this last kind are Nelder–Mead method, genetic algorithms, particle swarm optimization, simulated annealing, or ant colony optimization, see~\citet{arora_optimization_algorithms}.

\subsection{Separable nonlinear regression problems}
Some of the previous problems fall into the family of \textbf{separable nonlinear regression} problems. In this kind of problem, the nonlinear function $f$ in \eqref{eq:fullnonlinear} can be written as a linear combination of nonlinear functions. In particular, $f$ takes the expression
\begin{equation}
\label{eq:separablenls}
    f(\mathbf{x};\boldsymbol\theta) = \sum_{i = 1}^r a_i \phi_i(\mathbf{x};\mathbf{b}).
\end{equation}
With this formulation, the original set of parameters $\boldsymbol\theta$ has been split into two subsets: the linear parameters $\mathbf{a}=(a_1,\dots,a_r)$   and the nonlinear parameters $\mathbf{b}= (b_1,\ldots,b_s)$. Obviously $r+s = p$ and thus the number of nonlinear parameters to be determined is smaller than the total number of parameters. Therefore, the number of parameters to be found using nonlinear algorithms can be reduced. Separable nonlinear least squares methods are described, for example, in~\citet{golub_separable} and~\citet{golub_differentiation}.

In this work we present the package \pkg{nlstac}, based on the TAC algorithm described in~\citet{tac} for solving the separable nonlinear regression problem given in \eqref{eq:separablenls}.

\subsection{Related packages}
There are some widely used R packages that also deal with the problems \pkg{nlstac} is designed to solve. However, they rely on different algorithms which are often dependent on the choice of initial values. These packages are mainly \CRANpkg{nlsr} \citep{package_nlsr} and \pkg{minpack.lm}, both solving the problems with variants of the Levenberg-Marquardt algorithm, \CRANpkg{lbfgs} package \citep{package_lbfgs}, which provides an interface to L-BFGS and OWL-QN algorithms, or \CRANpkg{minqa} package \citep{package_minqa}, implementing derivative-free optimization algorithms. The algorithms used by this package fall into the aforementioned category of local algorithms.

Other R packages that use global algorithms, mostly of a stochastic nature, are, for example, \CRANpkg{DEoptimR} \citep{package_DEoptimR}, \CRANpkg{GenSA} \citep{package_GenSA}, \CRANpkg{GA} \citep{package_GA}, \CRANpkg{ABCoptim} \citep{package_ABCoptim} or \CRANpkg{pso} \citep{package_pso}.

There are also packages for fitting models in separable non-linear regression models, although these tend to be more specialized for specific problem domains. For example, the \CRANpkg{TIMP} package \citep{package_TIMP} is used for physics and chemistry problems whereas \CRANpkg{spant} package \citep{package_spant} deals with magnetic resonance spectroscopy problems. For partially separable nonlinear fitting we find \CRANpkg{psqn} package \citep{package_psqn}.

\section{The \pkg{nlstac} package}

The \pkg{nlstac} package was developed with two objectives: first, to implement the algorithm described in~\citet{tac} in functions that could be used for estimating separable nonlinear regression models, and second, to implement these functions with standard syntax such that they would be convenient for users familiar with other curve-fitting functions such as \code{lm}, \code{nls}, or \code{nlsLM} from the \pkg{minpack.lm} package.

The package consists of three units: a formula decomposer, a linear least squares solver, and a grid search unit. 

The workflow is depicted in Figure \ref{fig:workflow}: 

\begin{figure}[htp]
\includegraphics[width=\linewidth]{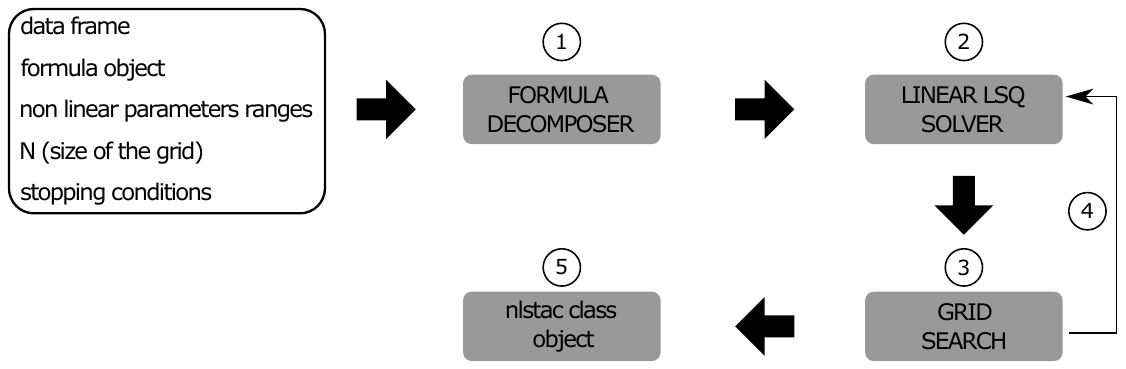}
\caption{Schematic workflow of the algorithm used in the \pkg{nlstac} package.}
\label{fig:workflow}
\end{figure}

\begin{enumerate}
\item From the dataset contained in a \code{dataframe},  an object of class \code{formula}, and a list of nonlinear parameters along with the initial ranges defined for each parameter, the formula decomposer, coded in internal function \code{get\_functions} determines the nonlinear functions, $\phi_i$, defined in \eqref{eq:separablenls}.
\item From the ranges (intervals) of each nonlinear parameter and the number of nodes, $N$, of each partition of such intervals, all possible combinations of nonlinear parameters are determined. For each combination, a linear least square problem is solved, obtaining thus a set of \textit{plausible parameters}. This is done by the  \code{get\_best\_parameters} internal function. 
Such a set of 'plausible parameters' is obtained in the following way: Let $b_1,\ldots,b_s$ be the nonlinear parameters in \eqref{eq:separablenls}. For each nonlinear parameter $b_i$ let $[c_i,d_i]$ denote the interval where to seek the estimation of $b_i$ and let  $c_i=b_i^1<b_i^2< \cdots < b_i^{N-1} < b_i^N =d_i$ denote the partition of such interval. Then, from these partitions, we construct a mesh in the rectangle $[c_1,d_1] \times \cdots \times [c_s,d_s]$. Next, for each node of the grid, we obtain the optimal linear parameters by solving a linear least square problem. The nodes of the grid, along with the optimal linear parameters for each node, constitute the set of \textit{plausible parameters}.
\item For each set of plausible parameters, the loss function (namely the sum of the squares of the residuals) is computed, and a grid search is performed to obtain the minimum value of the loss function. Let $(b_1^{m_1}, \ldots, b_s^{m_s})$ be the node of the grid in which the loss function minimizes.
\item  Stopping criteria are met when either the maximum number of iterations is reached or when the size of the partition of every interval $[c_i,d_i]$ is lower than the tolerance level. If stopping criteria are not met, the grid is refined: for each parameter, $b_i$, a new subinterval where to seek the estimation of such parameter is considered, $[b_i^{m_i-1},b_i^{m_i+1}]$ (note that if $m_i$ is either $1$ or $N$, the new subinterval will be $[c_i,b_i^2]$ or $[b_i^{N-1},d_i]$, respectively). The new grid will be established by repeating steps 2 to 4 until one stopping criterion is met.
\item When stopping criteria are met, the result is returned as an object of class \code{nlstac}. 
\end{enumerate}

As indicated in step 5, the output given by the \code{nls\_tac} function is an object of class \code{nlstac}. It is a list containing the following fields: 
\begin{itemize}
\item \code{resid}: The residuals.
\item \code{data}: The original data. 
\item \code{coefficients}: A named vector containing the values of the parameters.
\item \code{stdError}: A named vector with the standard error of the estimation of the coefficients.
\item \code{convInfo}: Convergence information. Namely, the number of iterations (\code{niter}), and the tolerance reached (\code{tolerance}).
\item \code{SSR}: The sum of the squares of the residuals obtained by the fit. 
\item \code{fitted}: A vector containing the fitted values.
\item \code{dataset}: A string with the name of the variable containing the data.
\item \code{formula}: The formula used in the call of the function.
\item \code{df}: The degrees of freedom
\item \code{sigma}: The standard deviation estimate.
\item \code{Rmat}: R matrix in the QR decomposition of the gradient matrix used for the computation of the standard errors of the coefficients.
\end{itemize}

The class \code{nlstac} has also some extraction methods similar to \code{lm}, \code{nls}, \code{glm}. For instance, the methods \code{summary.nlstac}, \code{predict.nlstac}, and \code{predict.summary.nlstac} produce identically formatted output as the \code{summary} functions for the \code{lm} and \code{nls} fits, as will be shown later. 

The \pkg{nlstac} package \citep{package_nlstac} is available in CRAN. The development version of the package can also be installed from the GitHub repository using the \code{install\_github} function from the remotes package \citep{package_remotes}: \code{remotes::install\_github("rbensua/nlstac")}.

\subsection{Arguments} \label{ss_arguments}
As was mentioned above, the inputs for the \code{nls\_tac} function are, at least, the fields \code{data}, \code{formula}, \code{tol}, \code{N} and \code{nlparam}.

The \code{data} field is a data frame containing  the data to be fitted; \code{tol} is the tolerance measured as the relative difference between the values of the parameters in two consecutive iterations; its  default value is $10^{-4}$; \code{N} is the number of divisions we make in each nonlinear parameter interval in each iteration (defaults to 10); \code{formula} is either an object of \code{formula} class or a character string that can be coerced into a \code{formula}, and it  must contain the pattern or formula which will be fitted, and \code{nlparam} is a list containing the nonlinear parameters as well as their initial ranges.

\subsection{Summary of nlstac class} \label{s_summary}
Information about the fit is stored in an \code{nlstac} class, and the function \code{summary} can display the most numerical relevant information about the fit.

Considering the example shown in Subsection \ref{ss_one_exp} (Example 1. Exponential Decay), once the analysis is done, the \code{summary} function shows us the information of the analysis as follows:
\begin{example}
  > summary(tacfit)

Formula: temp ~ a1 * exp(-k1 * time) + a2

Parameters:
       Estimate   Std. Error  t value   Pr(>|t|)    
k1 1.399458e-02 8.107657e-05 172.6095 < 2.22e-16 ***
a1 4.951112e+01 1.447617e-01 342.0182 < 2.22e-16 ***
a2 2.382372e+01 5.877739e-02 405.3212 < 2.22e-16 ***
---
Signif. codes:  0 ‘***’ 0.001 ‘**’ 0.01 ‘*’ 0.05 ‘.’ 0.1 ‘ ’ 1

Residual standard error: 0.1647017 on 180 degrees of freedom

Number of iterations to convergence:  13
Achieved convergence tolerance:  2.87864e-08
\end{example}

As can be seen, function \code{summary} gives us information about the formula used in the fitting, the estimated parameters along with some statistical information, the residual standard error, the number of iterations necessary to achieve convergence, and, finally, the tolerance value when convergence is achieved.

\section{Examples}
In this section, we present several examples to illustrate the use of the \pkg{nlstac} package in various scenarios. Each example highlights different aspects of \pkg{nlstac}'s behavior compared to another commonly used R function for these types of problems: \code{nls}, which is included in the \pkg{stats} package and utilizes the Gauss-Newton algorithm by default.

For each example, we explore two different initializations for \code{nls}. First, we initialize \code{nls} with a reasonable set of initial values. Note that the actual values of the parameters are unknown and no \textit{a priori} estimation of these values are available; therefore we denote as reasonable a set of values which are similar to the estimation obtained by TAC. This approach allows us to compare the fit achieved by \pkg{nlstac} with the widely used \code{nls} function. In the second initialization approach, we initialize \code{nls} with the estimation provided by \pkg{nlstac} for the same problem. This second approach serves as a starting point to ease the convergence of the algorithm used by \code{nls} and enables us to observe how effectively both algorithms can work in tandem. These examples shed light on the versatility and potential of the \pkg{nlstac} package in conjunction with the established \code{nls} function, providing valuable insights into their combined performance.

In Subsection \ref{ss_one_exp} and \ref{ss_two_exp} two examples with real data are presented. In both cases, the \pkg{nlstac} package obtains a solution, while the \code{nls} function does not converge with seemingly reasonable set of initial values. However, if the \code{nls} function is initialized with the output of \pkg{nlstac}, it successfully converges to a solution. In the first example, the fit remains the same, and in the second one, \code{nls} slightly improves the fit obtained by \pkg{nlstac}. These examples highlight the versatility of \pkg{nlstac}, which can be used either independently for fitting or to provide accurate initialization values for \code{nls} to converge effectively.

In Subsection \ref{ss_three_exp_phase} we present an example with simulated data. In this example, way beyond TAC proven convergence, we obtain a good fit to the model when running the \pkg{nlstac} package. However, if we initialize the \code{nls} function with seemingly reasonable set of initial values, it fails to converge. Nevertheless, by using the output of \pkg{nlstac} as initialization values for \code{nls}, we achieve an even better fit than what \pkg{nlstac} alone provides. This example also shows that \pkg{nlstac} can serve not only as a standalone fitting algorithm but also as a tool that enhances the performance of \code{nls} by providing reliable initialization values for improved fitting.

In Subsection \ref{ss_exp_sin}, we present another example with simulated data where the \pkg{nlstac} package accurately fits the given pattern. However, in this case, the \code{nls} function converges using both initialization approaches. Interestingly, when a non-optimal initialization is used, the fit obtained by \code{nls} is poor.

In Subsection \ref{ss_autorregressive}, we fit an exponential autoregressive model. This multi-variable example showcases the robustness of the TAC algorithm showing how the \pkg{nlstac} package can be utilized in wide range of scenarios. In this example, we generate three different datasets by making perturbations to a common underlying pattern. This example is quite interesting because for each dataset the behavior of \code{nls} differs.

Furthermore, in Subsection \ref{ss_autorregressive_real_data}, we utilize the \pkg{nlstac} package to fit real-world data to an autoregressive model.

In Subsection \ref{ss_function_explicit}, we illustrate how the \code{function} parameter in the \code{nls\_tac} function can be  utilized to explicitly provide the functions within the family of approximants. This feature proves useful in cases where the algorithm does not accurately recognize the pattern.

Lastly, we want to mention that all figures presented in this section have been generated using the \CRANpkg{ggplot2} package \citep{package_ggplot2}.

\newcounter{n}
\addtocounter{n}{1}
\subsection{Example \arabic{n}. Exponential Decay} \label{ss_one_exp}
Although we implement \pkg{nlstac} to fit data with virtually any nonlinear function, as mentioned in the introduction, the original purpose of the TAC algorithm was to fit exponential decays models such as \eqref{eq:pp}. The convergence of TAC for exponential decays patterns is proved using the max-norm as the approximation criterion. 

Patterns
\begin{equation*}
a_1 e^{-k_1t} + a_2, \; \text{where } a_1,a_2,k_1 \in \mathbb{R}, \; k_1>0,
\end{equation*}
presented in \eqref{eq:pp} are widely used to fit data coming from measuring the temperature of a body during a time interval, and their use for this propose is based on Newton’s law of cooling. Let us see an example of this use. 

Five parameters: \code{data}, \code{tol}, \code{N}, \code{formula} and \code{nlparam} need to be passed, as indicated in Subsection \ref{ss_arguments}. The first parameter, \code{data}, must be a 2-columns matrix containing data: instants and observations.

We intend to fit pattern \eqref{eq:pp} to dataset \textit{Coolingwater} from \CRANpkg{mosaicData} package \citep{package_mosaicData}. First, we define variable \code{data}.
\begin{example}
  data <- CoolingWater[40:222,]
\end{example}

Once \code{data} is loaded, we specify the tolerance, \code{tol}, or stopping criterion, and the number of divisions to be made in each step, \code{N}.
\begin{example}
  tol <- 1e-7
  N <- 10
\end{example}
We usually set the number of divisions to 10. However, if the search intervals for the nonlinear parameters are very wide or if we suspect that there may be many local minima, it might be advisable to increase the number of divisions to avoid converging to a sub-optimal local minimum. On the contrary, if we suspect that the computing time may be too high (for example, if the number of nonlinear parameters is large), it might be advisable to reduce the number of divisions.

Next, we specify the model to be used in the fitting, \code{form}, specifying the nonlinear parameters included in the model, \code{nlparam}, as well as the interval in which we assume they can be found. Please observe that the function does not require us to initialize the parameters whatsoever, we are just asked to provide a (wide) interval where to seek them. In this example, we have chosen the interval $[10^{-7},1]$ as the interval where $k_1$ must be sought.
\begin{example}
  form <- 'temp ~ a1*exp(-k1*time) + a2'
  nlparam <- list(k1 = c(1e-7,1))
\end{example}

Finally, we run the \code{nls\_tac} function to obtain the fit. 
\begin{example}
  tacfit <- nls_tac(formula = form, data = data,  nlparam = nlparam, N = N, tol = tol,
    parallel = FALSE)
\end{example}
Note that the input formula is either an R formula object or an object coercible to it. For example, in this case, variable form is a string that can be coerced to a formula object. Also note that we only need to specify the names and the initial intervals for the nonlinear parameters in the nlparam  input. Once the nonlinear parameters are given, the function \code{nls\_tac} will call the formula decomposer that will try to determine the rest of the elements of the formula ---i.e. the linear parameters and nonlinear functions described in equation \eqref{eq:separablenls}. Finally, note that \code{tacfit} is an object of class \code{nlstac} containing the following fields: \code{coefficients}, \code{stdError}, \code{convInfo}, \code{SSR}, \code{residuals}, \code{data} and \code{formula}.

So far we have only used one algorithm for the fitting: the TAC algorithm. A rightful question to be asked would be how much this fit resembles the one provided by the widely used nonlinear Least Squares (NLS) algorithm. 

For this purpose we will use the NLS algorithm, by means of \code{nls} function, to fit the same pattern to the same data. Since NLS requires initialization values, we initialize parameter $k_1$ as $0.1$, and parameters $a_1$ as $50$, and $a_2$ as $20$. Then we run \code{nls}.
\begin{example}
    nlsfit1 <- nls(formula = form, data = data, start = list(k1 = 0.1, a1 = 50, a2=20),
    control = nls.control(maxiter = 1000, tol = tol))
\end{example}
Although we have chosen reasonable initialization values, the algorithm did not converge. This shows a strength of TAC, which converges without the need for initialization values.

Besides the use of TAC and NLS as algorithms that can fit by themselves, they can also be used jointly. Since \pkg{nlstac} does not require initialization, just a set of intervals in which to seek the nonlinear parameters, \pkg{nlstac} can provide to \code{nls} good initialization values so that \code{nls} can successfully converge. The lack of dependence on initialization values of \pkg{nlstac} paired with the speed of \code{nls} and extended use among researchers, make them quite a good team.

When using \pkg{nlstac} and \code{nls} together, we use the coefficients obtained in the fit with \pkg{nlstac} to initialize and run \code{nls} for a second time.
\begin{example}
  nlsfit2 <- nls(formula = form, data = data, start = coef(tacfit), 
  control = nls.control(maxiter = 1000, tol = tol))
\end{example}
In this case, the \code{nls} did indeed converge, but the fit coincides with \pkg{nlstac}'s. This show that \code{nls} was not able to improve \pkg{nlstac} fit, even though it was initialized with its output.

We show the summary of \pkg{nlstac} and \code{nls} in Table  \ref{t:one_exponential_summaries}. For both methods, the residual standard error is 0.1647017 on 180 degrees of freedom. For \pkg{nlstac}, the number of iterations to convergence is 13 and the achieved convergence tolerance is 2.87864$\times 10^{-8}$. For \code{nls} when fitting with \pkg{nlstac} output as initialization, the number of iterations to convergence is 1 ---meaning the initial values were close to the optimal values; furthermore, the algorithm was unable to improve those values--- and the achieved convergence tolerance is 1.391929$\times 10^{-8}$.

\begin{table}[htb]
\centering
\begin{tabular}{ccccc}
\toprule
Parameter &  Estimate & Std. Error & t value & Pr(>|t|)\\   
\midrule
$k_1$ & 0.01399458 & 8.107657e-05 & 172.6095 & < 2.22$\times 10^{-16}$ ***\\
$a_1$ & 49.51112 & 0.1447617 & 342.0182 & < 2.22$\times 10^{-16}$***\\
$a_2$ & 23.82372 & 0.05877739 & 405.3212 & < 2.22$\times 10^{-16}$ ***\\
\bottomrule
\end{tabular} 
\caption{Example \arabic{n}. Summary of \pkg{nlstac} and \code{nls} for \textit{CoolingWater} dataset with the model given in \eqref{eq:pp}. Note that all outputs coincide for both methods.}
\label{t:one_exponential_summaries}
\end{table}

\begin{figure}[htb]
\centering
\includegraphics[width=0.8\textwidth]{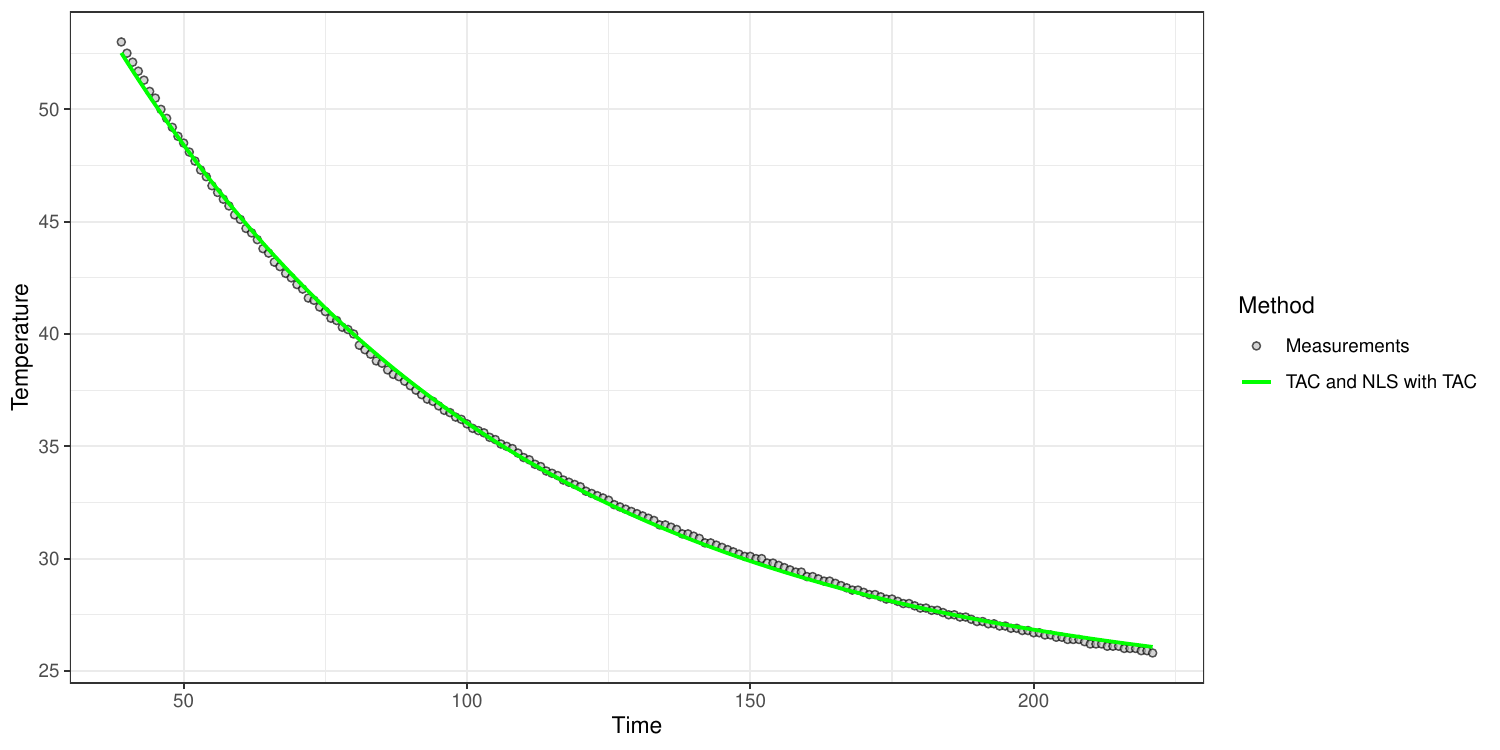}
\caption{Example \arabic{n}. Fitting an exponential decay for \textit{CoolingWater} dataset. The figure shows the original data (grey points) and the \pkg{nlstac} fit along with the \code{nls} fit (green line).}
\label{fig:one_exponential}
\end{figure}

Figure~\ref{fig:one_exponential} shows the data as gray dots. Both fits, the one provided by \pkg{nlstac} and the one provided by \code{nls} initialized using \pkg{nlstac}'s best approximation, are shown in green.  

For more examples of using TAC algorithm on real-world data, see section 4 of~\citet{tac} or subsection 5.1 of~\citet{tac_mdpi}.

\addtocounter{n}{1}
\subsection{Example \arabic{n}. Bi-exponential Decay}  \label{ss_two_exp}
In some approximation problems it is necessary to fit the sum of two exponential decays, as published in~\citet{beniteztoca2010}. While \pkg{nlstac} completely solves this problem, \code{nls}, given a bad initial values, does not converge. However, if we initialize \code{nls} with \pkg{nlstac} output, \code{nls} fit the data in a very similar way that \pkg{nlstac} does.

In this example, we intend to fit a function such as
\begin{equation}\label{eq:2_exp}
f(t)=a_1 e^{-k_1t} + a_2 e^{-k_2t}  + a_3, \; \text{where } a_1,a_2,a_3,k_1,k_2 \in \mathbb{R}, k_1,k_2>0.
\end{equation}
Models of the form such as \eqref{eq:2_exp} were used in~\citet{beniteztoca2010} in the fitting of data produced in indentation experiments carried out by scanning probe microscopes (e.g., Atomic Force Microscopes) in studies of viscoelastic mechanical properties of soft matter.

We intend to fit pattern \eqref{eq:2_exp} to \textit{Indometh} data from the \pkg{datasets}  package \citep{R}. We define parameter \code{data} and specify the tolerance, \code{tol}, and the number of divisions made in each step, \code{N}:
\begin{verbatim}
  data <- Indometh[Indometh$Subject == 3, ]
  tol <- 1e-7
  N <- 10
\end{verbatim}

We set the model to be used in the fitting, \code{form}, specifying the nonlinear parameters included in the model, \code{nlparam}, as well as the interval in which we assume they can be found. Finally, we apply the \code{nls\_tac} function to get the fit. 
\begin{example}
  form <- 'conc ~ a1*exp(-k1*time) + a2*exp(-k2*time) + a3'
  nlparam <- list(k1 = c(1e-7,10), k2 = c(1e-7,10))
  tacfit <- nls_tac(formula = form, data = data,  nlparam = nlparam, N = N, tol = tol,
    parallel = FALSE)
\end{example}
In a similar way as indicated in Example \ref{ss_one_exp}, we run \code{nls} initializing every parameter, that is, $k_1$, $k_2$, $a_1$, $a_2$ and $a_3$, as 1. Later we use the coefficients obtained with \pkg{nlstac} to initialize and run \code{nls} for a second time.
\begin{example}
  nlsfit1 <- nls(formula = form, data = data, 
    start = list(k1 = 1,k2 = 1, a1 = 1, a2 = 1, a3 = 1), 
    control = nls.control(maxiter = 1000, tol = tol))
  nlsfit2 <- nls(formula = form, data = data, start = coef(tacfit), 
    control = nls.control(maxiter = 1000, tol = tol))
\end{example}
While running this last piece of code we encountered an error because of bad initial values when using a vector of ones to initialize \code{nls}. We get a gradient error and \code{nls} does not converge. However, if parameters are initialized using \pkg{nlstac} output, \code{nls} does converge. 

We show the summaries of \pkg{nlstac} and \code{nls} in Table \ref{t:two_exponential_summary_TAC} and Table \ref{t:two_exponentials_summary_NLS_best_approximation}, respectively.

\begin{table}[htb]
\centering
\begin{tabular}{ccccc}
\toprule
Parameter &  Estimate & Std. Error & t value & Pr(>|t|)\\    
\midrule     
$k_1$ & 0.94813244 & 0.14656766 & 6.46891 & 0.00064760 ***\\
$k_2$ & 9.75308642 &  5.32150841 & 1.83277 & 0.11654040 \\
$a_1$ & 2.11675784  & 0.29303012 & 7.22369 & 0.00035686 ***\\
$a_2$ & 11.09789500 & 12.99197213 & 0.85421 & 0.42577278\\
$a_3$ & 0.08168448  & 0.03165979 & 2.58007 & 0.04176551 * \\
\bottomrule
\end{tabular} 
\caption{Example \arabic{n}. Summary of \pkg{nlstac} for data of subject 3 in \textit{Indometh} dataset with model given in \eqref{eq:2_exp}. Residual standard error: 0.05351802 on 6 d.o.f. Number of iterations to convergence: 12. Achieved convergence tolerance:  3.824026$\times 10^{-8}$.}
\label{t:two_exponential_summary_TAC}
\end{table}

\begin{table}[htb]
\centering
\begin{tabular}{ccccc}
\toprule
Parameter &  Estimate & Std. Error & t value & Pr(>|t|)\\    
\midrule    
$k_1$  & 0.89971458 & 0.14914067 & 6.03266 & 0.00093743 ***\\
$k_2$ & 7.96454599 & 3.19653381 & 2.49162 & 0.04705893 *  \\
$a_1$ & 2.00446255 & 0.29689160 & 6.75150 & 0.00051489 ***\\
$a_2$ & 7.63334977 & 4.94487377 & 1.54369 & 0.17361052    \\
$a_3$ & 0.07663298 & 0.03262943 & 2.34858 & 0.05716893 .  \\
\bottomrule
\end{tabular} 
\caption{Example \arabic{n}. Summary of \code{nls} for data of subject 3 in \textit{Indometh} dataset with model given in \eqref{eq:2_exp}. Residual standard error: 0.0527844 on 6 d.o.f. Number of iterations to convergence: 8. Achieved convergence tolerance:  3.646981$\times 10^{-8}$.}
\label{t:two_exponentials_summary_NLS_best_approximation}
\end{table}

\begin{figure}[htb]
\centering
\includegraphics[width=0.8\textwidth]{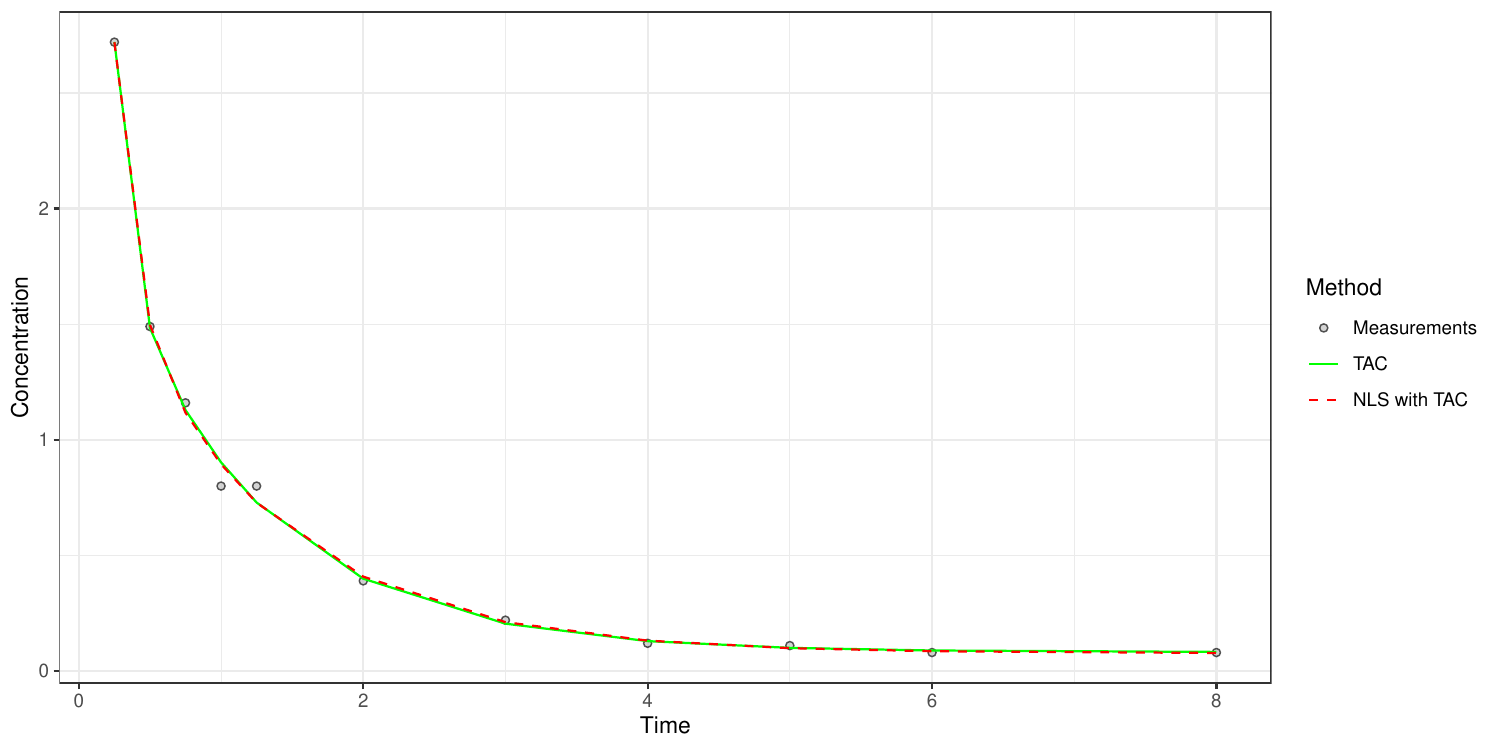}
\caption{Example \arabic{n}. Fit of a bi-exponential decay for data of subject 3 in \textit{Indometh} dataset. The figure shows the original data (grey points), the \pkg{nlstac} fit (green line) and the \code{nls} fit (red line)}
\label{fig:two_exponentials}
\end{figure}

Figure~\ref{fig:two_exponentials} shows the data as gray dots; \pkg{nlstac} fit is shown in green and, in dashed red, \code{nls} fit is shown. 

Although it is now hidden from the user, this implementation used to show some warnings related to a deficiency in a matrix rank. The explanation is that we are assuming both parameters $k_1$ and $k_2$ are contained within an interval $[10^{-7},10]$, so when we make two equal partitions of the same interval, one for each parameter, at some point both values will be the same: nodes of $[10^{-7},10] \times [10^{-7},10]$ where $k_1=k_2$ and therefore we do not have a bi-exponential decay but only an exponential decay. Since these values are used in the resolution of a linear equation system that does not have a unique solution, we obtain a warning because such a unique solution can not be found. The algorithm takes care of that problem removing these indeterminate parameter sets, as well as the warnings.

Also note that, for some choice of models, permutations of parameters may give the same results. For example, using the model above, for each pair of parameters $(k_1,k_2)$ there will be another pair of parameters $(k_2,k_1)$ which will offer the same fit. However, that is not a problem for us, other than for an increase in the time of execution of the code. Another similar scenario is when adjusting two sinusoidal waves and two exponential decays. If we wanted to avoid making the same calculations multiple times, we would have to change the code, forcing the user to specify which functions are non-identifiable for permutations of parameters, so we would get a more time-efficient code at the cost of simplicity. However we have chosen simplicity over time efficiency.

Another bi-exponential decay example with real data can be found in section 5 of~\citet{tac}.

\addtocounter{n}{1}
\subsection{Example \arabic{n}. Exponential decays with phase displacement} \label{ss_three_exp_phase}
In this example, \pkg{nlstac} converges and \code{nls}, even when reasonable initialization values are given, does not. However, if we use \pkg{nlstac} output as an initialization for \code{nls},  \code{nls} not only converges but improves \pkg{nlstac} fitting.

Here we get to see two advantages of TAC algorithm. Firstly, \pkg{nlstac} converges. Secondly, when given its approximation for \code{nls} initialization, \code{nls} improves upon this fit. That shows us, again, the two ways of using \pkg{nlstac}: directly estimating models or providing initial values.

We intend to fit a function such as
\begin{equation} \label{eq:three_exponentials}
f(t)=a_1 e^{-b_1(t-d_1)^2} + a_2 e^{-b_2(t-d_2)^2} + a_3 e^{-b_3(t-d_3)^2}  + a_4,
\end{equation}
where $a_1,a_2,a_3,a_4,b_1,b_2,b_3,d_1,d_2,d_3 \in \mathbb{R}$. As indicated in Subsection \ref{ss_arguments} we need to pass five parameters: \code{data}, \code{tol}, \code{N}, \code{formula} and \code{nlparam}.

We create \code{data} and determine the tolerance, \code{tol}, and the number of divisions we make in each step, \code{N}.
\begin{example}
  set.seed(12345)
  x <- seq(from = 0, to = 20, length.out = 65)
  y <- 2*exp(-10*(x-0.5)^2) + 3*exp(-1*(x-2)^2) + 4*exp(-0.1*(x-5)^2) + 1 + .05*rnorm(65)
  data <- data.frame(x,y)
  tol <- 1e-5
  N <- 6
\end{example}

We set the number of divisions to 6 because otherwise the computation time for \pkg{nlstac} is too high due to the increased number of nonlinear parameters. This is one of the downsides of \pkg{nlstac}, the increase of time processing is too high when several nonlinear parameters need to be fitted.

Later, we specify the model to be fitted, \code{form}, specifying the nonlinear parameters included in the model, \code{nlparam}, as well as the intervals in which we assume they can be found. 
\begin{example}
  form <- 'y ~ a2*exp(-b1*(x-d1)^2) + a3*exp(-b2*(x-d2)^2)+ a4*exp(-b3*(x-d3)^2)+ a1'
  nlparam <- list(b1 = c(7.7,15), b2 = c(0,5.1), b3 = c(1e-4,1.1),
                d1 = c(1e-2,1.5), d2 = c(0.1,4), d3 = c(0.11,11))  
\end{example}
Finally, we apply the \code{nls\_tac} function to get the fit. 
\begin{example}
    tacfit <- nls_tac(formula = form, data = data,  nlparam = nlparam, N = N, tol = tol,
    parallel = FALSE)
\end{example}

As indicated in previous examples, we will run \code{nls}. First, we will initialize it with the following values: $a_1 = 0.5$, $a_2 = 1$, $a_3 = 3$, $a_4 = 5$, $b_1 = 10$, $b_2 = 0.5$, $b_3 = 0.1$, $d_1 = 0$, $d_2 = 1$ and $d_3 = 1$. Later, we will provide \pkg{nlstac} output as \code{nls} initialization
\begin{example}
    nlsfit1 <- nls(formula = form, data = dat, 
      start = list(a1 = 0.5, a2 = 1, a3 = 3, a4 = 5, b1 = 10, b2 = 0.5, b3 = 0.1, 
      d1 = 0, d2 = 1, d3 = 1), control = nls.control(maxiter = 1000, tol = tol))
    nlsfit2 <- nls(formula = form, data = dat, start = coef(tacfit), 
      control = nls.control(maxiter = 1000, tol = tol))
\end{example}
As commented before, \code{nls} does not converge with the first initialization but does converge with \pkg{nlstac} initialization. 

We show the results of both implementations in Table \ref{t:three_exponentials}.

\begin{table}[htb]
\centering
\begin{tabular}{cccccc}
\toprule
Method & $a_1$ & $a_2$ & $a_3$ & $a_4$  &  $b_1$ \\
\midrule
\pkg{nlstac} & 1.00896226 & -1.5798685 & 3.38784985 & 3.95157736 & 7.7\\
\code{nls} & 1.00522131  & -2.00224548 & 3.79346203 & 3.98814647 & 3.86837135\\
\midrule
Method & $b_2$ & $b_3$  &  $d_1$ & $d_2$ & $d_3$  \\
\midrule
\pkg{nlstac} & 0.53116510 & 0.10359427 & 1.21392 & 1.64537609 & 5.08417547\\
\code{nls} & 0.61731929 & 0.09950229 & 1.26058001 & 1.57407878  & 5.01184086 \\
\bottomrule
\end{tabular} 
\caption{Example \arabic{n}. Parameter estimates corresponding to \pkg{nlstac} and \code{nls} fit for dataset considered in example \arabic{n} with the model given in \eqref{eq:three_exponentials}. Values have been rounded to the eighth decimal place. Please recall that \code{nls} initialized without the estimate from \pkg{nlstac} does not converge.}
\label{t:three_exponentials}
\end{table}

We show the summaries of \pkg{nlstac} and \code{nls} in Table \ref{t:three_exponentials_summary_TAC} and Table \ref{t:three_exponentials_summary_NLS_best_approximation}, respectively.

\begin{table}[htb]
\centering
\begin{tabular}{ccccc}
\toprule
Parameter &  Estimate & Std. Error & t value & Pr(>|t|)\\   \midrule     
$a_1$ & 1.008962  & 0.026108 & 38.646 & < 2$\times 10^{-16}$ ***\\
$a_2$ & -1.579868 &  0.200716 & -7.871 & 1.41$\times 10^{-10}$ ***\\
$a_3$ & 3.387850  & 0.189397 & 17.888 & < 2$\times 10^{-16}$ ***\\
$a_4$ & 3.951577  & 0.060719 & 65.080 & < 2$\times 10^{-16}$ ***\\
$b_1$ & 7.700000  & 2.148880 &  3.583 & 0.00072 ***\\
$b_2$ & 0.520107  & 0.030797 & 16.888 & < 2$\times 10^{-16}$ ***\\
$b_3$ & 0.103594  & 0.058464 &  1.772 & 0.08195 .  \\
$d_1$ & 1.213920  & 0.048245 & 25.162 & < 2$\times 10^{-16}$ ***\\
$d_2$ & 1.645376  & 0.007954 & 206.864 & < 2$\times 10^{-16}$ ***\\
$d_3$ & 5.084175  & 0.099682 & 51.004 & < 2$\times 10^{-16}$ ***\\
\bottomrule
\end{tabular} 
\caption{Example \arabic{n}. Summary of \pkg{nlstac} for dataset considered in example \arabic{n} with the model given in \eqref{eq:three_exponentials}. Residual standard error: 0.141 on 55 d.o.f. Number of iterations to convergence: 14. Achieved convergence tolerance:  9.871$\times 10^{-6}$.}
\label{t:three_exponentials_summary_TAC}
\end{table}

\begin{table}[htb]
\centering
\begin{tabular}{ccccc}
\toprule
Parameter &  Estimate & Std. Error & t value & Pr(>|t|)\\   \midrule    
$a_1$ & 1.005221 &  0.024285 & 41.392 & < 2$\times 10^{-16}$ ***\\
$a_2$ &-2.002245 &  0.350988 & -5.705 &4.80$\times 10^{-7}$ ***\\
$a_3$ & 3.793462 &  0.311814 & 12.166 & < 2$\times 10^{-16}$ ***\\
$a_4$ & 3.988146 &  0.054770 & 72.816 & < 2$\times 10^{-16}$ ***\\
$b_1$ & 3.868371 &  1.061791 &  3.643 &0.000597 ***\\
$b_2$ & 0.617319 &  0.078514 &  7.863 &1.46$\times 10^{-10}$ ***\\
$b_3$ & 0.099502 &  0.006762 & 14.715 & < 2$\times 10^{-16}$ ***\\
$d_1$ & 1.260580 &  0.033273 & 37.886 & < 2$\times 10^{-16}$ ***\\
$d_2$ & 1.574079 &  0.047734 & 32.976 & < 2$\times 10^{-16}$ ***\\
$d_3$ & 5.011841 &  0.087899 & 57.018 & < 2$\times 10^{-16}$ ***\\
\bottomrule
\end{tabular} 
\caption{Example \arabic{n}. Summary of \code{nls} for dataset considered in example \arabic{n} with the model given in \eqref{eq:three_exponentials}. Residual standard error: 0.1307 on 55 d.o.f. Number of iterations to convergence: 16. Achieved convergence tolerance:  2.8$\times 10^{-6}$.}
\label{t:three_exponentials_summary_NLS_best_approximation}
\end{table}

\begin{figure}[htb]
\centering
\includegraphics[width=0.8\textwidth]{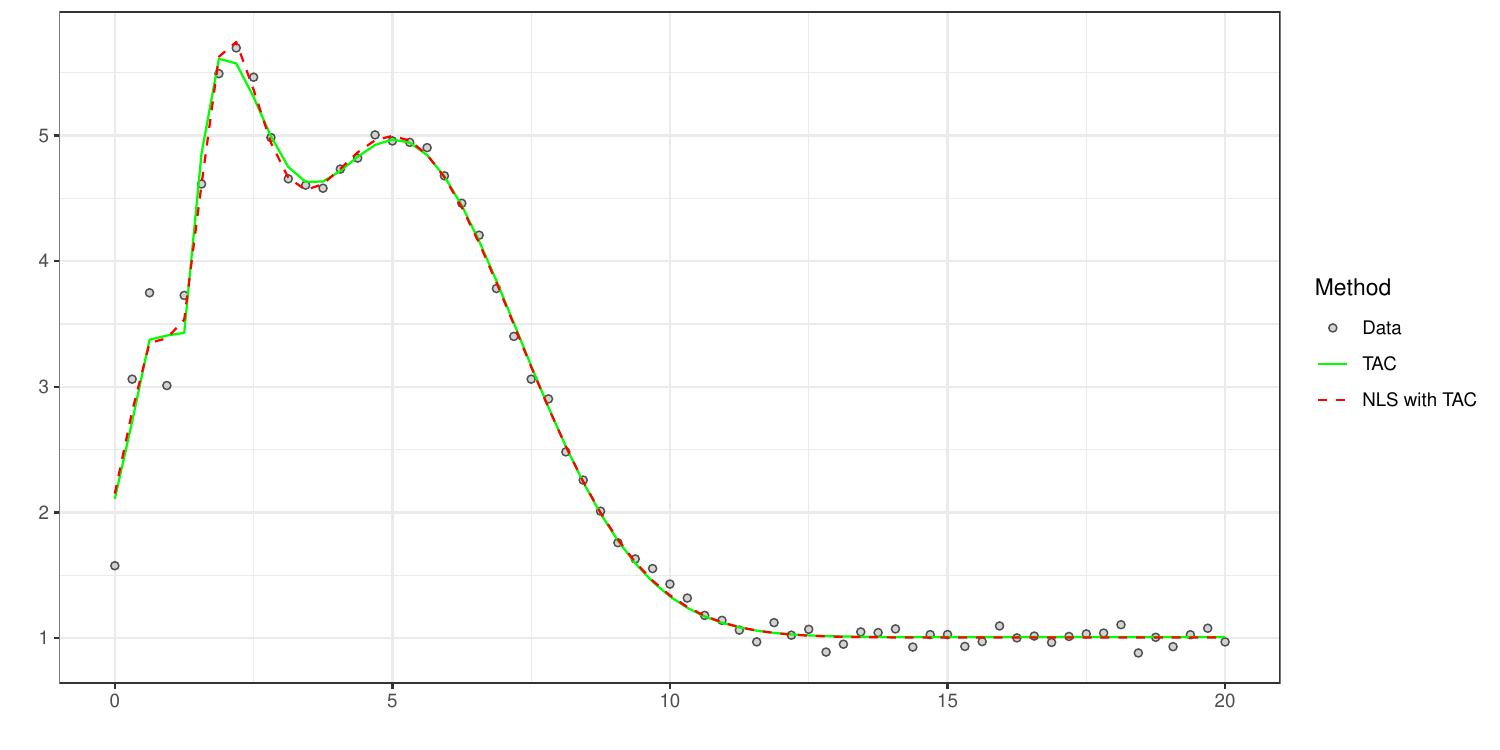}
\caption{Example \arabic{n}. The combined trend of three exponential decays with phase displacement. The figure shows the original data (grey points), the \pkg{nlstac} fit (green line) and the \code{nls} fit using \pkg{nlstac}'s best approximation (red line).}
\label{fig:three_exponentials}
\end{figure}

Figure~\ref{fig:three_exponentials} shows the data as gray dots. In green, \pkg{nlstac} fit is shown. Dashed red line shows \code{nls} fit initialized with the parameters of \pkg{nlstac}'s best approximation.

This example is particularly significant since \pkg{nlstac} is outperformed by \code{nls}, both in time (\pkg{nlstac}'s computing time is significantly higher) and precision (compare the residual standard error value for both methods in Tables \ref{t:three_exponentials_summary_TAC} and \ref{t:three_exponentials_summary_NLS_best_approximation}), although \code{nls} function needs to be initialized with \pkg{nlstac} best approximation to be able to converge.

\addtocounter{n}{1}
\subsection{Example \arabic{n}. Exponential decay mixed with a sinusoidal signal}\label{ss_exp_sin}
In this example, where we mix an exponential decay with a sinusoidal signal, we obtain a good fit with \pkg{nlstac} and a similar fit with \code{nls} when we provide \pkg{nlstac} output as initialization values. However, if we initialize \code{nls} with a bad initialization values, we get a poor fit: \code{nls}' fit identifies the exponential decay quite properly but fails to identify the sinusoidal signal.

We intend to fit a function such as
\begin{equation} \label{eq:pattern_exp_sin}
f(t)=a_1 e^{-k_1 t} + a_2 \sin(b_1 t) + a_3, \; \text{where } a_1,a_2,a_3,k_1,b_1 \in \mathbb{R}, k_1,b_1>0.
\end{equation}
As indicated in \ref{ss_arguments} we need to pass five parameters: \code{data}, \code{tol}, \code{N}, \code{formula} and \code{nlparam}.

We create \code{data} and determine the tolerance, \code{tol}, or stopping criterion, and the number of divisions to be made in each step, \code{N}.
\begin{example}
  set.seed(12345)
  x <- seq(from = 0, to = 10, length.out = 500)
  y <- 3*exp(-0.85*x) + 1.5*sin(2*x) + 1 + rnorm(length(x), mean = 0, sd = 0.3)
  data <- data.frame(x,y)
  tol <- 1e-7
  N <- 10
\end{example}

Later we set the model to be used in the fitting, \code{form}, specifying the nonlinear parameters included in the model, \code{nlparam}, as well as the intervals in which we assume they can be found.
\begin{example}
 form <- 'y ~ a1*exp(-k1*x) + a2*sin(b1*x) + a3'
 nlparam <- list(k1 = c(0.1,1), b1 = c(1.1,5))
\end{example}
Finally, we apply the \code{nls\_tac} function to adjust the data. 
\begin{example}
  tacfit <- nls_tac(formula = form, data = data,  nlparam = nlparam, N = N, tol = tol,
    parallel = FALSE)
\end{example}

As in previous examples, we compare the \pkg{nlstac} and \code{nls} output. For the first comparison we will run \code{nls} initializing every parameter, that is, $k_1$, $b_1$, $a_1$, $a_2$ and $a_3$, as 1, and for the second comparison, we will use  \pkg{nlstac} output to initialize and run \code{nls}.
\begin{example}
    nlsfit1 <- nls(formula = form, data = data, start = list(k1 = 1, b1 = 1, a1 = 1,
      a2 = 1, a3 = 1) , control = nls.control(maxiter = 1000))
    nlsfit2 <- nls(formula = form, data = data, start = coef(tacfit), control = 
      nls.control(maxiter = 1000))
\end{example}

We present the results obtained in Table \ref{t:one_exponential_one_sinusodial}. Without looking at any graphic, it is quite evident from Table \ref{t:one_exponential_one_sinusodial} alone that the last fit is different from the two others.

\begin{table}[htb]
\centering
\begin{tabular}{cccccc}
\toprule
Method & $k_1$ & $b_1$ & $a_1$ & $a_2$ & $a_3$ \\ 
\midrule
\pkg{nlstac} &  0.81234568  & 1.99851628  & 3.01166130  & 1.51095645 & 1.00766609 \\ 
\midrule
\code{nls}  &  &  &  & &  \\ 
(with \pkg{nlstac} &  0.81904149 & 1.99847422 & 3.01996411 & 1.51073313 & 1.00969794\\ 
initialization) &  &  &  & &  \\ 
\midrule
\code{nls}  &  &  &  & &  \\ 
(without \pkg{nlstac} &  0.82435634 & 0.83083729 & 4.49199629 & -0.23904801 &  0.91892869 \\
initialization) &  &  &  & &  \\  
\bottomrule
\end{tabular} 
\caption{Example \arabic{n}. Parameters corresponding to \pkg{nlstac} and \code{nls} fits for dataset considered in example \arabic{n} with the model given in \eqref{eq:pattern_exp_sin}. Values have been rounded off to the eighth decimal place.}
\label{t:one_exponential_one_sinusodial}
\end{table}

Summary of \pkg{nlstac} is shown in Table \ref{t:one_exponential_one_sinusodial_summary_TAC} and Table \ref{t:one_exponential_one_sinusodial_summary_NLS_best_approximation} shows the summary of \code{nls} initialized with the best approximation obtained with \pkg{nlstac}. Finally, summary of \code{nls} initialized with a vector of ones appears in Table \ref{t:one_exponential_one_sinusodial_summary_NLS_bad_approximation}.

\begin{table}[htb]
\centering
\begin{tabular}{ccccc}
\toprule
Parameter &  Estimate & Std. Error & t value & Pr(>|t|)\\    
\midrule     
$k_1$ & 0.812346  & 0.035145 &  23.11  & <2$\times 10^{-16}$ ***\\
$b_1$ & 1.998516 &  0.002132 & 937.35 & <2$\times 10^{-16}$ ***\\
$a_1$ & 3.011661 &  0.077098 &  39.06 &  <2$\times 10^{-16}$ ***\\
$a_2$ & 1.510956 &  0.019734 &  76.57  & <2$\times 10^{-16}$ ***\\
$a_3$ & 1.007666  & 0.018885 &  53.36 &  <2$\times 10^{-16}$ ***\\
\bottomrule
\end{tabular} 
\caption{Example \arabic{n}. Summary of \pkg{nlstac} for dataset considered in example \arabic{n} with the model given in \eqref{eq:pattern_exp_sin}. Residual standard error: 0.2974 on 495 d.o.f. Number of iterations to convergence: 11. Achieved convergence tolerance:  3.184$\times 10^{-8}$.}
\label{t:one_exponential_one_sinusodial_summary_TAC}
\end{table}

\begin{table}[htb]
\centering
\begin{tabular}{ccccc}
\toprule
Parameter &  Estimate & Std. Error & t value & Pr(>|t|)\\    
\midrule   
$k_1$ & 0.819041 &  0.035388 &  23.14  & <2$\times 10^{-16}$ ***\\
$b_1$ & 1.998474 &  0.002132 & 937.31  & <2$\times 10^{-16}$ ***\\
$a_1$ & 3.019964 &  0.077382  & 39.03 &  <2$\times 10^{-16}$ ***\\
$a_2$ & 1.510733 &  0.019733  & 76.56 &  <2$\times 10^{-16}$ ***\\
$a_3$ & 1.009698 &  0.018813 &  53.67 &  <2$\times 10^{-16}$ ***\\
\bottomrule
\end{tabular} 
\caption{Example \arabic{n}. Summary of \code{nls} initialized with \pkg{nlstac}'s best approximation for dataset considered in example \arabic{n} with the model given in \eqref{eq:pattern_exp_sin}. Residual standard error: 0.2974 on 495 d.o.f. Number of iterations to convergence: 4. Achieved convergence tolerance:  5.597$\times 10^{-9}$.}
\label{t:one_exponential_one_sinusodial_summary_NLS_best_approximation}
\end{table}

\begin{table}[htb]
\centering
\begin{tabular}{ccccc}
\toprule
Parameter &  Estimate & Std. Error & t value & Pr(>|t|)\\    
\midrule      
$k_1$ & 0.82436  &  0.10069 &  8.187 & 2.3$\times 10^{-15}$ ***\\
$b_1$  & 0.83084 &   0.06042 & 13.751 & < 2$\times 10^{-16}$ ***\\
$a_1$ & 4.49200  &  0.26879 & 16.712 & < 2$\times 10^{-16}$ ***\\
$a_2$ & -0.23905 &   0.07680 & -3.112 & 0.00196 ** \\
$a_3$ & 0.91893 &   0.07455 & 12.326 & < 2$\times 10^{-16}$ ***\\
\bottomrule 
\end{tabular} 
\caption{Example \arabic{n}. Summary of \code{nls} initialized with a vector of ones for dataset considered in example \arabic{n} with the model given in \eqref{eq:pattern_exp_sin}. Residual standard error: 1.056 on 495 d.o.f. Number of iterations to convergence: 23. Achieved convergence tolerance:  7.587$\times 10^{-8}$.}
\label{t:one_exponential_one_sinusodial_summary_NLS_bad_approximation}
\end{table}

\begin{figure}[htb]
\centering
\includegraphics[width=0.8\textwidth]{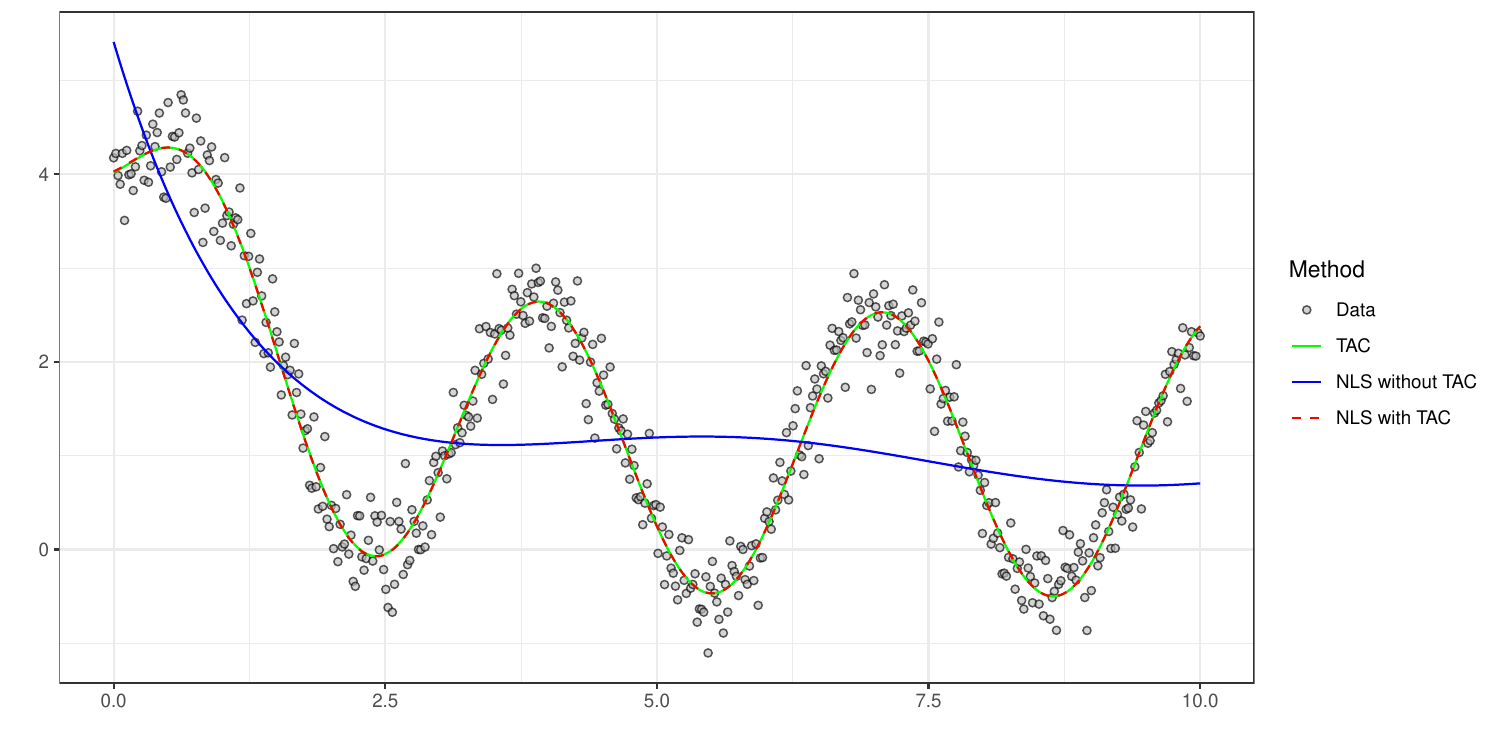}
\caption{Example \arabic{n}. Fit of an exponential decay mixed with a sinusoidal signal for dataset considered in example \arabic{n} with the model given in \eqref{eq:pattern_exp_sin}. The figure shows the original data (grey points), the \pkg{nlstac} fit (green line), the \code{nls} fit initialized with \pkg{nlstac} output (red line) and the  \code{nls} fit initialized with a vector of ones (blue line).}
\label{fig:one_exponential_one_sinusodial}
\end{figure}

Figure~\ref{fig:one_exponential_one_sinusodial} shows the data as gray dots. Green line represents the \pkg{nlstac} fit and dashed red line represents \code{nls} fit. Blue line represents \code{nls} fit initialized with a vector of ones. It is clear that the last fit is not accurate. It seems that the nonlinear least squares algorithm has managed to fit the exponential part of the pattern but it seems to have missed the sinusoidal part.

This example shows a situation where \pkg{nlstac} works perfectly, as well as \code{nls} if initialized correctly. However, if the user does not provide good initialization values, the nonlinear least squares algorithm might fail to obtain a good fit since it may get stuck in a local minimum.
 
\addtocounter{n}{1}
\subsection{Example \arabic{n}. Exponential autoregressive model: a multi-variable approach (p-variable)}\label{ss_autorregressive}
Nonlinear time series models are used in a wide range of fields. In this example, we deal with an especially relevant nonlinear time series model: the exponential autoregressive model. Given a time series $\{ x_1,x_2, x_3,\ldots \}$, the exponential autoregressive model is defined as
\begin{equation*}
 x_t= \left[ \sum_{i=1}^p \left(a_i + b_i e^{-cx_{t-1}^2}\right)x_{t-i} \right] + \varepsilon_t,
\end{equation*}
where $\varepsilon_t$ are independent and identically distributed random variables and independent with $x_i$, $p$ denotes the system degree, $t \in \mathbb{N}$, $t>p$, and the parameters to be estimated from observations are $c$, $a_i$ and $b_i$ (for $i=1,\ldots,p$). This model can be found in, for example,~\citet{example_xu_modeling_a_nonlinear} or~\citet{example_chen_generalized_exponential_autoregressive}.

Some generalizations for the exponential autoregressive model have been made, and in this example, we will deal with a generalization of Teräsvirta’s extended model that can be found in~\citet[equation (10)]{example_chen_generalized_exponential_autoregressive} and we present here:
\begin{equation}\label{eq:ExpAR_Terasvirta_generalization}
 x_t= a_0 + \left[ \sum_{i=1}^p \left(a_i + b_i e^{-c (x_{t-d} \; - \; z_i)^2}\right)x_{t-i} \right] + \varepsilon_t,
\end{equation}
where $z_i$ (for $i=1,\ldots,p$) are scalar parameters and $d \in \mathbb{Z}$.  

We would like to point out that the convergence of TAC algorithm has not been established for this type of problem. Further, this example is substantially different from the above examples since every observation depends on the previous ones. This model can not be described by a function of just one real variable. Instead, a vector of $p$ real variables needs to be used. This approach can be developed considering a function from $(\mathbb{R}^{n-p})^p$ into $\mathbb{R}^{n-p}$. Therefore we transform a one-dimensional problem into a $p$-dimensional one. Let us explain this process. Let $x=(x_1,\ldots,x_n)$ denote the observations and let us define $p$ variables $v1,\ldots,vp \in \mathbb{R}^{n-p}$, being $vi=(x_{p-i+1}, \ldots, x_{n-i})$ for $i=1,\ldots,p$, which will allow us to redefine Equation \eqref{eq:ExpAR_Terasvirta_generalization} in terms of these new variables:
\begin{equation}\label{eq:ExpAR_Terasvirta_generalization_with_p_variables}
 x_{t+p}= a_0 + \left[ \sum_{i=1}^p \left(a_i + b_i e^{-c ((vd)_t \; - \; z_i)^2}\right)(vi)_t \right] + \varepsilon_t, \text{ with } t=1,\ldots,n-p,
\end{equation}
where $vd$ is fixed with $d$ such that $1 \leq d \leq p$ and $(vi)_t$ denote the $t$-th component of variable $vi$.

For this example, first, we are going to fix one exponential autoregressive time series. Then we are going to generate three different datasets: Dataset 1, Dataset 2, and Dataset 3. Each of these datasets is generated by setting three different seeds (Seed 1, Seed 2, and Seed 3, respectively) in order to add a random perturbation to the terms of the previously fixed time series.

The aim of this example is to fit model \eqref{eq:ExpAR_Terasvirta_generalization_with_p_variables} with $p=2$, $d=2$, and $n=100$ for each dataset. We start defining a vector, \code{seed}, containing the three seeds. We define the tolerance level, \code{tol}, and the parameters for the time series as well as initialize the time series with the first two terms. We also define the pattern to be fitted.
\begin{example}
  seed <- c('12','123','1234')
  tol <- 1e-5
  a0 <- -1.45
  a1 <- 1.66
  b1 <- -0.47
  a2 <- 0.543
  b2 <- -0.82
  c <- 1.27
  z1 <- 2.53
  z2 <- 3.85
  x <- numeric(100)
  x[1] <- 2.7
  x[2] <- 3.12
  form <- 'y ~ a0+ a1*v1 + b1*v1*exp(-c*(v2-z1)^2) + a2*v2 + b2*v2*exp(-c*(v2-z2)^2)'
\end{example}

We intend to run \pkg{nlstac} in parallel using package \CRANpkg{doParallel} \citep{package_doParallel}. We set up the parallelization:
\begin{example}
  no_cores <- detectCores() - 1
  cl <- makeCluster(no_cores)
  registerDoParallel(no_cores)
\end{example}

We are going to define a loop that will iterate three times, one for each dataset. In each iteration we set a seed and generate its corresponding dataset: seeds '12', '123', and '1234' will generate Dataset 1, Dataset 2, and Dataset 3, respectively. Then we transform the problem into a two-dimensional one by defining variables \code{y}, \code{v1}, and \code{v2} as previously described. We create the dataframe, run \pkg{nlstac} in parallel and run \code{nls} with two different initializations: first initialized with $c = 1$, $z_1 = 2.25$, $z_2 = 4$, $a_0 = -1$, $a_1 = 1$, $b_1 = -1$, $a_2 = 1$, $b_2 = -1$ and then initialized with \pkg{nlstac} output. The function \code{tryCatch} is used in order to keep the loop running in the event that \code{nls} does not converge for some choice of initial parameters. 
\begin{example}
  for (j in 1:3) {
    set.seed(seed[j])
    for (i in 3:100){
      x[i] <- a0 + (a1+b1*exp(-c*(x[i-2]-z1)^2))*x[i-1] + 
        (a2+b2*exp(-c*(x[i-2]-z2)^2))*x[i-2] + rnorm(1, mean=0, sd=0.1)}
    y <- x[3:100]
    v1 <- x[2:99]
    v2 <- x[1:98]
    data <- data.frame(v1 = v1, v2 = v2, y = y)
    tacfit <- nls_tac(form, data = data, 
                      nlparam = list(c = c(0,2), z1 = c(1,3.5), z2 = c(3,5)),
                      N = 10, tol=tol, parallel = TRUE)
    tryCatch(nlsfit1 <- nls(formula = form, data = data, start = list(c = 1, z1 = 2.25, 
    	z2 = 4, a0 = -1, a1 = 1, b1 = -1, a2 = 1, b2 = -1), 
    	control = nls.control(maxiter = 1000, tol = tol, minFactor = 1e-5)), 
    	error = function(e) {nlsfit1 <<- NULL})
    tryCatch(nlsfit2 <- nls(formula = form, data = data, start = coef(tacfit), 
    	control = nls.control(maxiter = 1000, tol = tol, minFactor = 1e-5)), 
    	error = function(e) {nlsfit2 <<- NULL}) }
\end{example}

Finally we stop the parallelization:
\begin{example}
  stopImplicitCluster()
\end{example}

For Data 1, pattern \eqref{eq:ExpAR_Terasvirta_generalization_with_p_variables} has successfully been fitted with all three methods. For Data 2, \pkg{nlstac} did converge but \code{nls} did not converge for either of the two initializations. Finally, for Data 3, \code{nls} only converged if initialized with \pkg{nlstac} output, and in this case \code{nls} fit slightly improves \pkg{nlstac}'s.

We present the results for all three datasets in Table \ref{t:expAR_several_seeds}. 

\begin{widetable}[htb]
\centering
\begin{tabular}{*{8}{c}} 
\toprule
  &   &  \multicolumn{3}{c}{Dataset 1}  & Dataset 2  & \multicolumn{2}{c}{Dataset 3}\\
  \cmidrule(lr){3-5}
  \cmidrule(lr){6-6}
  \cmidrule(lr){7-8}
   &      &  & \code{nls} & \code{nls} &  &  &  \code{nls}  \\    
 \multicolumn{2}{c}{Parameter}     & \pkg{nlstac} & (without  & (with \pkg{nlstac} & \pkg{nlstac} &  \pkg{nlstac} & (with \pkg{nlstac}  \\  
      &      &  & \pkg{nlstac} init.) & init.) &  &   &  init.) \\  
\midrule     
\multirow{2}{*}{$c$} & Estimate & 1.038256  &  1.019213  &    1.019233  &  0.246914  &  1.728395  & 5.441692  \\
  & Std. Error & 0.487366  &  0.477747  &  0.477748  &  3.816922  &  2.300992 &  4.654574 \\
\midrule 
\multirow{2}{*}{$z_1$} & Estimate & 2.388889  &  2.436606  &   2.436604  &  3.211423  &  1.438958 &   2.261739  \\
  & Std. Error & 0.521781  &  0.501099  &  0.501096 &   0.082645  &  1.707290  &  0.287022 \\
\midrule 
\multirow{2}{*}{$z_2$} & Estimate & 3.855393  &  3.881704   &   3.881695 &   3.006097  &  3.843975  &  3.628858   \\
  & Std. Error & 0.251267  &  0.231838  &  0.231837  &  2.613194  &  0.296154  &   0.123689 \\
\midrule 
\multirow{2}{*}{$a_0$} & Estimate & -1.929269  & -2.691265 &  -2.691134  &  5.616140  &  -8.359206  &  -3.401572   \\
  & Std. Error & 7.038082  &  6.687379  &   6.687395  &  45.308458  & 9.373772  &   5.347638 \\
\midrule 
\multirow{2}{*}{$a_1$} & Estimate & 1.778348  &  1.823764  &    1.823752   & -3.721976  & 1.312138 &   1.329813 \\
  & Std. Error & 0.409774  &  0.437640   &    0.437630  &   77.352777   &  0.075125 &   0.067164  \\
\midrule 
\multirow{2}{*}{$b_1$} & Estimate & -0.813623  &  -0.832519 &  -0.832507  &  5.169974  &  2.888960 &  0.327253   \\
  & Std. Error & 0.360567  &  0.380349  &    0.380342  &  77.339916  &   11.657008  &   0.542132 \\
\midrule 
\multirow{2}{*}{$a_2$} & Estimate & 1.055334  &  1.349577  &    1.349517  &  1.794767  &  2.716872  & 0.851547      \\
  & Std. Error & 2.698594  &   2.549761 &    2.549766  &   40.429455  &  3.641374   &  1.854555 \\
\midrule 
\multirow{2}{*}{$b_2$} & Estimate & -1.272387  &  -1.392871 &  -1.392838  &  -4.026071  &  -0.873735  &  -0.310493 \\
  & Std. Error & 1.115143  &  1.090259  &   1.090251  &  55.467877   &  1.228695  &  0.348730 \\
\midrule 
\multicolumn{2}{c}{SSR} & 0.701835 & 0.701740 & 0.701740 & 0.690280 & 0.87582 & 0.856081 \\
\bottomrule
\end{tabular} 
\caption{Example \arabic{n}. Summary of all three methods (\pkg{nlstac}, \code{nls} without \pkg{nlstac} initialization, \code{nls} with \pkg{nlstac} initialization) for all three datasets considered in example \arabic{n} with the model given in \eqref{eq:ExpAR_Terasvirta_generalization}. Missing methods for Dataset 2 and Dataset 3 are the result of the non-convergence of such methods.}
\label{t:expAR_several_seeds}
\end{widetable}

\begin{figure}
\centering
\includegraphics[width=\linewidth]{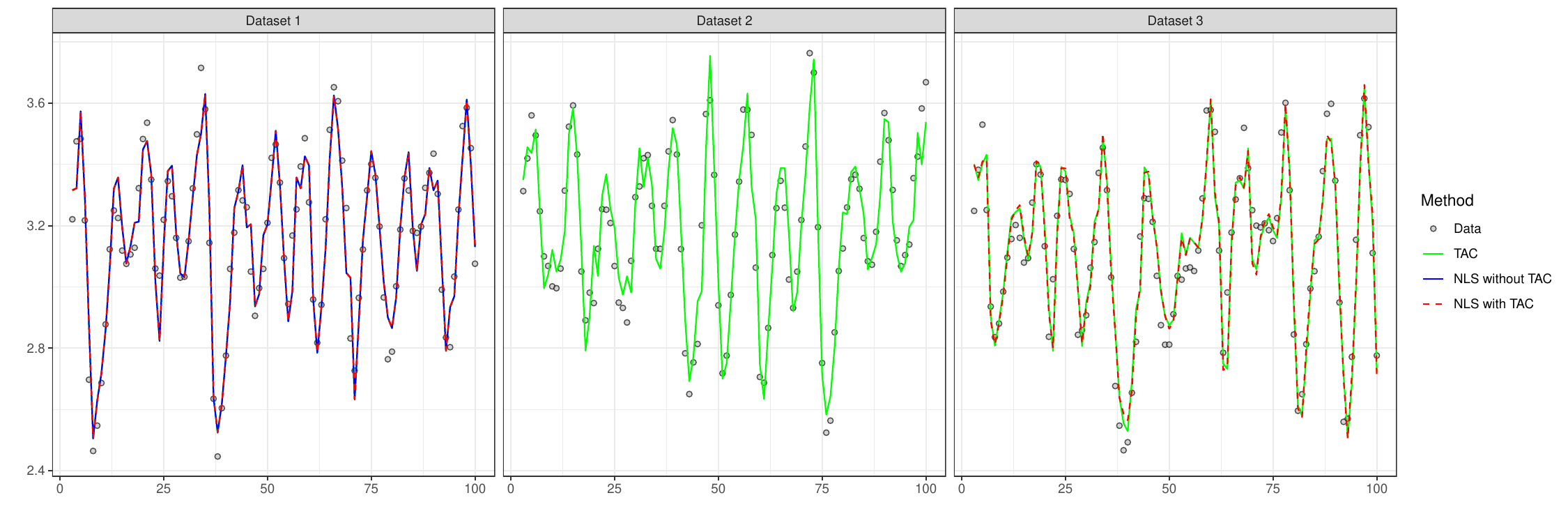}
\caption{Example \arabic{n}. Fits for three exponential autorregresive model datasets with model given in \eqref{eq:ExpAR_Terasvirta_generalization} as determined at the beginning of Example \arabic{n}. The figure shows the original data (grey points), the \pkg{nlstac} fit (green line), the \code{nls} fit initialized with \pkg{nlstac} output (red line) and the \code{nls} fit initialized with \pkg{nlstac} output (blue line). Note that, in Dataset 1, blue line overlaps the green one.}
\label{fig:expAR}
\end{figure}

Figure~\ref{fig:expAR} shows the fitting for those three datasets. 

\addtocounter{n}{1}
\subsection{Example \arabic{n}. Exponential autoregressive model: a multi-variable approach (p-variable) with real data.}\label{ss_autorregressive_real_data}
In this example, we intend to fit model \eqref{eq:ExpAR_Terasvirta_generalization} to returns on daily closing prices data from Financial Times Stock Exchange (FTSE) using \pkg{nlstac} package. More precisely, if vector $x=(x_1,\ldots,x_n)$ denotes the daily closing prices data, we are going to fit the returns, that is, vector $(\frac{x_2 - x_1}{x_1}, \ldots, \frac{x_i - x_{i-1}}{x_{i-1}}, \ldots , \frac{x_n - x_{n-1}}{x_{n-1}} )$. This data was obtained by \textit{EuStockMarkets} dataset which is accessible from \pkg{datasets} package. As in the previous example, we are going to consider $p=2$ and $d=2$.

First, we read in the data and specify the tolerance level:
\begin{example}
  x <- EuStockMarkets[,4]
  x <- diff(x)/x[-length(x)]
  tol <- 1e-7
\end{example}

Then we transform the problem into a two-dimensional one:
\begin{example}
  y <- x[3:length(x)]
  v1 <- x[2:(length(x)-1)]
  v2 <- x[1:(length(x)-2)]
  data <- data.frame(v1 = v1, v2 = v2, y = y)
  form <- 'y ~ a0+ a1*v1 + b1*v1*exp(-c*(v2-z1)^2) + a2*v2 + b2*v2*exp(-c*(v2-z2)^2)'
\end{example}

Finally, we make use of package \pkg{doParallel} to run \pkg{nlstac} in parallel:
\begin{example}
  no_cores <- detectCores() - 1
  cl <- makeCluster(no_cores)
  registerDoParallel(no_cores)
  tacfit <- nls_tac(form, data = data, 
                  nlparam = list(c = c(1e-7,5), z1 = c(1e-7,5), z2 = c(1e-7,5)),
                  N = 15, tol=tol, parallel = TRUE)
  stopImplicitCluster()
\end{example}

Results are summarized in Table \ref{t:expAR_summary_TAC_real_data} and a plot with both the data and the fit obtained by  \code{nls\_tac} function is depicted in Figure~\ref{fig:expAR} 

\begin{table}[htb]
\centering
\begin{tabular}{ccccc}
\toprule
Parameter &  Estimate & Std. Error & t value & Pr(>|t|)\\    
\midrule     
$c$  &  2.857144e-01  &  4.246881e+02  &  0.00067  &  0.99946  \\
$z_1$  & 3.668653e-03  &  4.806933e-03  &  0.76320  &  0.44544  \\
$z_2$ &  7.164164e-02  &  2.281940e-01  &  0.31395  &  0.75359  \\
$a_0$ &  7.571841e-05  &  2.219138e-04  &  0.34121  &  0.73299  \\
$a_1$ & -9.276088e+02  &  1.378917e+06  &  -0.00067  &  0.99946  \\
$b_1$ &  9.277265e+02  &  1.378917e+06  &  0.00067  &  0.99946  \\
$a_2$ & -1.374292e+02  &  2.038046e+05  &  -0.00067  &  0.99946  \\
$b_2$ &  1.376215e+02  &  2.038049e+05  &  0.00068  &  0.99946  \\
\bottomrule
\end{tabular} 
\caption{Example \arabic{n}. Summary of \pkg{nlstac} for returns from FTSE in \textit{EuStockMarkets} dataset with the model given in \eqref{eq:ExpAR_Terasvirta_generalization}. Residual standard error: 0.00791658 on 1849 d.o.f. Number of iterations to convergence: 11. Achieved convergence tolerance: 1.94289$\times 10^{-16}$.}
\label{t:expAR_summary_TAC_real_data}
\end{table}

\begin{figure}[htb]
\centering
\includegraphics[width=0.7\textwidth]{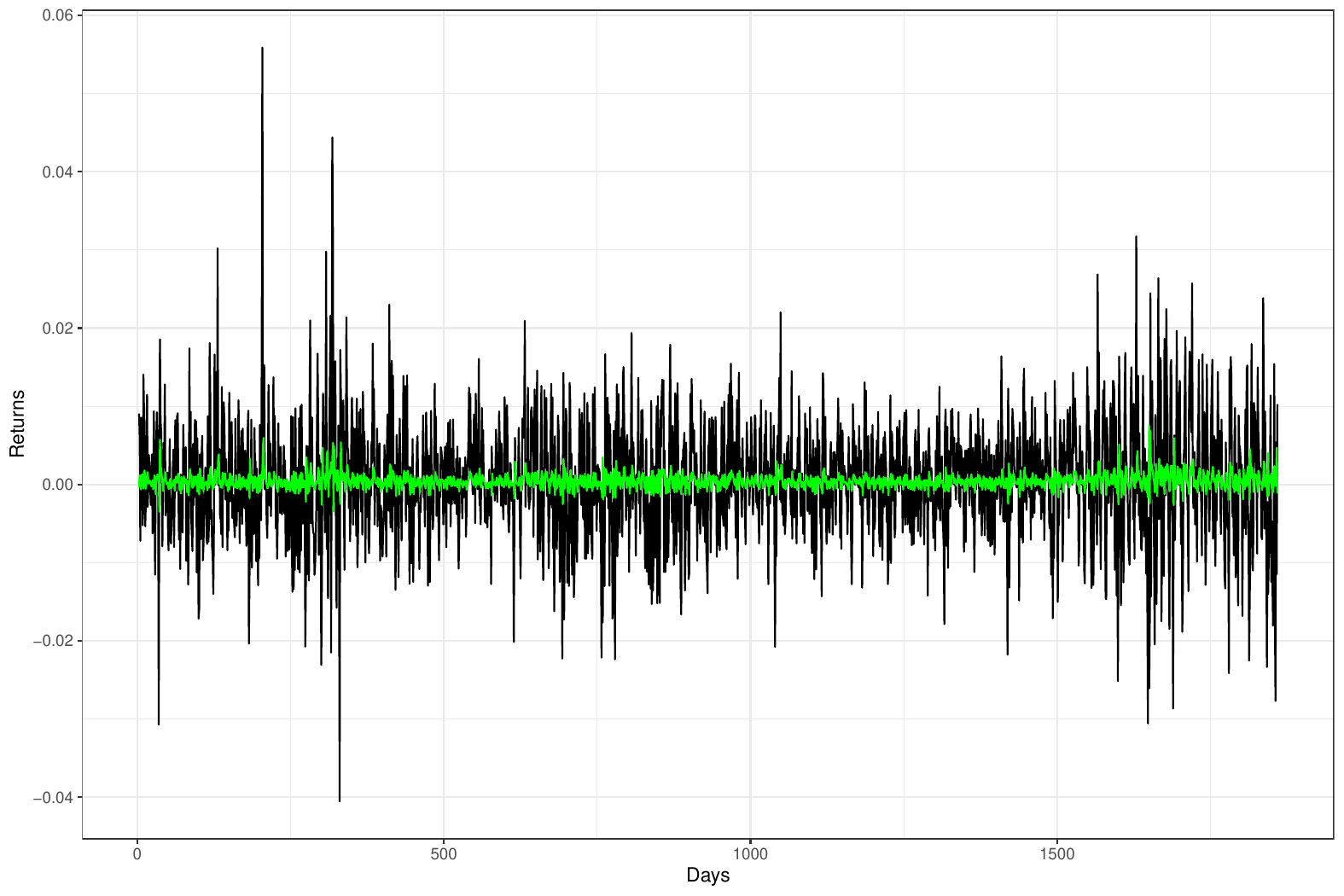}
\caption{Example \arabic{n}. Fit of an exponential autorregresive model for returns from FTSE in \textit{EuStockMarkets} dataset with model given in \eqref{eq:ExpAR_Terasvirta_generalization}. The figure shows the original data (black lines) and the \pkg{nlstac} fit (green line).}
\label{fig:expAR_real_data}
\end{figure}

Another example of an exponential autoregressive model with real data can be consulted in subsection 6.3 of~\citet{tac_mdpi}.

\addtocounter{n}{1}
\subsection{Example \arabic{n}. Explicitly providing the functions of the pattern}\label{ss_function_explicit}
Package \pkg{nlstac} relies on the \code{formula} infrastructure of the \code{R} language to determine the number and expressions of the nonlinear functions $\phi_i$ in \eqref{eq:separablenls}. This is done internally by the function \code{get\_functions}. However, in some cases, the user may want to explicitly state which are the nonlinear functions that define the separable nonlinear problem (e.g. the formula decomposer fails to automatically identify those functions). In that case, there is an optional parameter in the \code{nls\_tac} function, named \code{functions}, which is an array of character strings defining the nonlinear functions. In practical terms, we have not found an example in which we needed to manually specify the functions defining the model, but the optional parameter is available nonetheless.

Next, just for illustration purposes, we present here the example shown in Subsection \ref{ss_exp_sin} adding the \code{functions} parameter to the function \code{nls\_tac} so we explicitly provide the functions:

\begin{example}
  set.seed(12345)
  x <- seq(from = 0, to = 10, length.out = 500)
  y <- 3*exp(-0.85*x) + 1.5*sin(2*x) + 1 + rnorm(length(x), mean = 0, sd = 0.3)
  data <- data.frame(x,y)
  tol <- 1e-7
  N <- 10
  form <- 'y ~ a1*exp(-k1*x) + a2*sin(b1*x) + a3'
  nlparam <- list(k1 = c(0.1,1), b1 = c(1.1,5))
  tacfit <- nls_tac(formula = form, data = data, 
    functions=c('exp(-k1*x)','sin(b1*x)','1'),  nlparam = nlparam, N = N, tol = tol, 
    parallel = FALSE)
\end{example}

\section{Code parallelization}

The basic idea of the TAC algorithm is to find the optimal values for the linear parameters (by means of the linear least-square method) for each combination of the nonlinear parameters. Therefore, for every such combination of nonlinear parameters we have to solve a completely independent optimization problem, and thus this algorithm can take advantage of parallelization. 

The \pkg{nlstac} package implements a parallelization of this stage of the algorithm in the \code{nls\_tac} function. Setting the option \code{parallel=TRUE}, the function makes use of the \code{\%dopar\%} and \code{foreach} functions of the \CRANpkg{foreach} \citep{package_foreach} package and the infrastructure provided by the \pkg{parallel} \citep{R} and \pkg{doParallel} packages.

One might think that parallelization always speeds up the algorithm, but in reality, initializing and stopping the cluster requires a certain amount of time. Therefore, in some cases it may be convenient to parallelize and in others it might not be worth it. 

As was mentioned in the Introduction, the TAC algorithm, as all grid-search algorithms, scales poorly with the dimension of the problem (i.e. the number of nonlinear parameters). However, even for low-dimensional problems, the speed of the algorithm depends on the number of subdivisions of each parameter search interval (i.e. the width of the grid), which is defined by the parameter \code{N} in the \code{nls\_tac} function. Previous affirmations rely on the fact that, for each iteration, the number of plausible nonlinear parameters happens to be $N^p$, representing $N$, for each parameter, the number of nodes belonging to the partition of the interval where the parameter is assumed to be in and $p$ the number of nonlinear parameters. Note that the number of plausible nonlinear parameters depends on $N$ and exponentially increases with the number of nonlinear parameters $p$.

\begin{figure}[htb]
    \centering
    \includegraphics[width=0.9\linewidth]{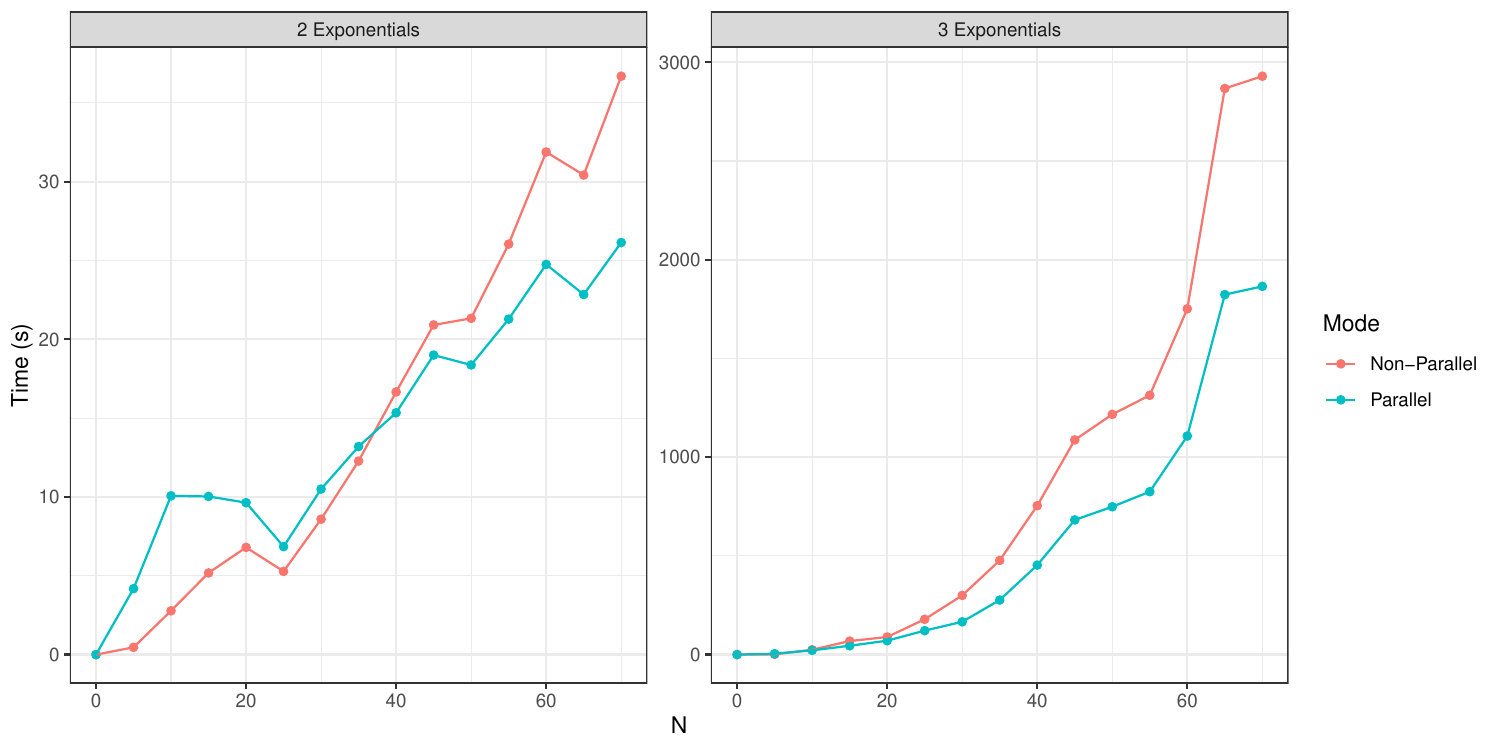}
    \caption{Comparison between the parallel and the non-parallel implementations of the \code{nls\_tac} function for the fitting of two (left) and three (right) exponential decays, for different values of the number of subdivisions, \code{N}. Note how the time of processing is increased around a hundred times when changing from two nonlinear parameters to three nonlinear parameters.}
    \label{fig:bench23}
\end{figure}

As an illustration of the convenience of using the \code{parallel = TRUE} option, Figure~\ref{fig:bench23} depicts a comparison of non-parallel and parallel modes of the \code{nls\_tac} function for two standard separable nonlinear least square problems. Namely, two exponential decays and three exponential decays. That is: 
\begin{equation*}
y = a_0 + \sum_{k=1}^n a_k e^{-b_k x},\qquad n = 2,3.
\end{equation*}
As could be expected, we can see how  as the number of nonlinear parameters increases, the computation time rises exponentially. Also, it shows that only for very small problems (e.g. two nonlinear parameters) the parallelization is not worth it in some cases (\code{N} up to 35). To run this simulation we had to make use of \CRANpkg{dplyr} package \citep{package_dplyr}.

\section{Conclusions}
Many popular packages for nonlinear function estimation depend heavily on the choice of starting values. This package, however, implements an algorithm that needs no initialization and can handle a wide variety of approximation problems.

Our goal has been to create a package for nonlinear regression using the TAC algorithm and to show how this algorithm can work either by itself or when combined with other nonlinear estimation algorithms. 

Processing times on problems with a large number of nonlinear parameters can be a problem. In those cases, it might be advisable to consider the use of a gradient-based algorithm. In future versions, the implemented grid search could be refined to reduce those processing times.

Despite this possible drawback, we strongly believe that this package will be found useful by researchers in nonlinear regression problems.


\bibliography{torvisco-benitez-arias-cabello.bib}

\address{J. A. F. Torvisco. Facultad de Ciencias. Universidad de Extremadura. Avda de Elvas s/n. 06006 Badajoz. Spain. ORCiD: 0000-0001-8373-3477. \email{jfernandck@alumnos.unex.es}}

\address{R. Ben\'{\i}tez. Departmento de Matemáticas para la Econom\'{\i}a y la Empresa. Universidad de Valencia. Avda Tarongers s/n, 46022 Valencia. ORCiD: 0000-0002-9443-0209. 
\email{rabesua@uv.es}}

\address{M. R. Arias. Facultad de Ciencias. Universidad de Extremadura. Avda de Elvas s/n. 06006 Badajoz. Spain. ORCiD: 0000-0002-4885-4270. \email{arias@unex.es}}
  
\address{J. Cabello~Sánchez. Facultad de Ciencias. Universidad de Extremadura. Avda de Elvas s/n. 06006 Badajoz. Spain. ORCiD: 0000-0003-2687-6193. \email{coco@unex.es}}

\end{article}

\end{document}